\documentclass[12pt]{amsart}
\usepackage{amssymb,amscd}
\usepackage{verbatim}
\usepackage{graphicx}

\newtheorem*{Whitney towers}{Theorem~\ref{Whitney towers}}
\newtheorem*{h-towers}{Theorems ~\ref{half} \& \ref{$(n)$-solvable}}

\newtheorem*{surgery curves}{Theorem~\ref{surgery curves}}
\newtheorem*{cg=0}{Theorem~\ref{vanish}}

\newtheorem{thm}{Theorem}[section]

\newtheorem{cor}[thm]{Corollary}

\newtheorem{cla}[thm]{Claim}
\theoremstyle{definition}
\newtheorem{defn}[thm]{Definition}

\newtheorem{prob}[thm]{Problem}

\numberwithin{equation}{section}
\numberwithin{figure}{section}

\newcommand{\x}{\times}
\newcommand{\np}{\newpage}
\newcommand{\Z}{\mathbb{Z}}
\newcommand{\N}{\mathbb{N}}

\newcommand{\R}{\mathbb{R}}

\newcommand{\f}{\noindent}

\begin{document}

\title{A new invariant and decompositions of manifolds}
\author{Eiji Ogasa\\}
\thanks{
Keywords: decomposition of manifolds, 
a new invariant $\nu(M)$, 
boundary union.   
\newline MSC2000 57N10, 57N13, 57N15. 
}
\date{}

\begin{abstract} 
We introduce a new topological invariant $\in\N\cup\{0\}$ 
of compact manifolds with boundaries 
associated with a kind of decomposition of them. 
Let $M$ and $N$ be 
$m$-dimensional compact connected manifolds with boundaries. 
Let $M\cup_{\partial}N$ be a boundary union of $M$ and $N$. 
Let $\nu(M)$ be the new invariant of $M$.  
Then we have \newline
$0\leqq\nu(M\cup_{\partial}N)\leqq \mathrm{max}\{\nu(M),\nu(N)\}.$
\end{abstract} 
\maketitle

\section{A problem}\label{introduction}

\f In order to state our problem we prepare a definition. 
We work in the smooth category. 

\begin{defn}\label{bdrysum}
Let $M$ (resp. $N$) be an $m$-dimensional 
smooth connected compact manifold with boundary. 
Let $\partial M=\amalg_{i=1}^{i=\alpha} X_i$ and 
$\partial N=\amalg_{j=1}^{j=\beta} Y_j$, 
where $\amalg$ denotes a disjoint union and $\alpha, \beta\in\N$. 
Suppose that each $X_i$ (resp. $Y_j$) is a connected closed manifold. 
A {\it boundary union}  $M\cup_{\partial}N$ 
is  an $m$-manifold with boundary 
which is a union $M\cup N$ with the following properties:

\begin{enumerate}
\item 
There are closed manifolds 
$X_{\sigma_1}\amalg...\amalg X_{\sigma_\mu}\subset\partial M$ and  
$Y_{\tau_1}\amalg...\amalg Y_{\tau_\mu}\subset\partial N$ 
which are diffeomorphic. 
Here, 
$\sigma_*$ are different each other,  
$\tau_\sharp$ are different each other, 
$1\leqq\sigma_*\leqq\alpha$,  
$1\leqq\tau_\sharp\leqq\beta$, 
$\mu\in\N$, 
$\mu\leqq\alpha$,  and  
$\mu\leqq\beta$.

\item 
We make $M\cup_{\partial}N$ 
by identifying  
$X_{\sigma_1}\amalg...\amalg X_{\sigma_\mu}$ 
with  
$Y_{\tau_1}\amalg...\amalg Y_{\tau_\mu}$ 
by a diffeomorphism map 
$X_{\sigma_1}\amalg...\amalg X_{\sigma_\mu}\to Y_{\tau_1}\amalg...\amalg Y_{\tau_\mu}$. 

\noindent
(Hence, of course, $M\cap N
=X_{\sigma_1}\amalg...\amalg X_{\sigma_\mu}
= Y_{\tau_1}\amalg...\amalg Y_{\tau_\mu}$.)  

\end{enumerate}

Let $\rho$ be an integer $\geqq2$.   
Suppose that  
a boundary union $L'$ of 
$\rho$ `manifolds with boundaries',   
$L_1,...,L_\rho$,  
is defined. 
Then a boundary union of 
`manifolds with boundaries', 
$L_{\rho+1}$ and $L'$,  
is said to 
be a {\it boundary union of 
$(\rho+1)$ `manifolds with boundaries',} 
$L_1,...,L_{\rho+1}$, 
be denoted by 
\hskip1mm
$\begin{matrix}
\text{{\tiny $\rho+1$}}\\
\text{\large $\cup_{\partial}$}\\
\text{{\tiny $i=1$}}\\
\end{matrix}$
\hskip1mm
$L_i$.
%
%
%
We say that $M$ is a {\it boundary union of one connected `manifold with boundary',} $M$.  
We say that  
the disjoint union $M\amalg N$ is a boundary union of $M$ and $N$. 
\end{defn}


\smallbreak
\noindent
{\bf Note.}  
Not all unions are boundary unions.   
Let  $B^3$ denote the 3-ball.  
Let  $D^2$ denote the 2-disc. 
We can regard  the solid torus as a union $B^3\cup B^3$   
such that $B^3\cap B^3$ is a disjoint union $D^2\amalg D^2$.  
This union $B^3\cup B^3$ is not a boundary union because  $B^3\cap B^3$ is 
not a disjoint union of closed manifolds. 

An example of boundary unions is the following: 
Let $T^2$ denote the torus. 
Let $A$ denote the annulus. 
Let $D^\circ$ denote the open 2-disc.    
Let $S^1$ denote the circle. 
We can regard $T^2-D^\circ$ as a union of 
$A\cup (A-D^\circ)$ such that 
$A\cap (A-D^\circ)$ is a disjoin union $S^1\amalg S^1$.  
This union $A\cup (A-D^\circ)$ is a boundary union. 

Any connected sum is a boundary union. However, the converse is not true. 
A \newline
Heegaard splitting of a closed oriented 3-manifold gives a boundary union of two copies of a handle body, which consists of a single 3-dimensional 0-handle and 3-dimensional 1-handles. 
However, it is not a connected sum if the handle body is not the 3-ball.

\bigbreak
We state our problem. 

\begin{prob}\label{ourproblem}    
Let $m$ be a nonnegative integer. 
Is there a finite set $\mathcal S$ of 
(oriented) compact connected 
$m$-manifolds with boundaries 
with  the following property ($\star$)?

\noindent
$(\star)$  
For an arbitrary closed connected $m$-manifold $M$, 
there are  `$m$-manifolds with boundaries'  $Z_1,...,Z_\nu$, where $\nu\in\N$, 
such that $M$ is a boundary union of  $Z_1,...,Z_\nu$ 
and that $Z_i\in\mathcal S$. 
Note that $Z_i$ may be same as $Z_j$ if $i\neq j$.


\end{prob}

Of course, we can make a problem in the $\partial M\neq\phi$ case    
if we impose an appropriate condition on $\partial M$. 

\smallbreak
It is trivial that the answer is affirmative if $m\leqq2$. 
String theory uses the fact 
that the $m=2$ case has the affirmative answer,  
discussing the world sheet (see  \cite{GSW, Polchinski} etc). 

At least, to the author, a motivation of this paper is the following: 
In QFT,  each Feynman diagram is made by the given fundamental parts. 
In string theory a world sheet (so to say, 2-dimensional Feynman diagram)
is decomposed into a finite number of 2-manifolds as stated above. 
In  $M$-theory  we may need high dimensional Feynman diagrams 
(see  \cite{Polchinski},     
P. 607, 608 of \cite{Vafa} and so on)  
although,  of course, there is an obstruction written in P60 of I of \cite{GSW}.
Considering high dimensional Feynman diagrams, 
we need to research a kind of decomposition of manifolds,   
for example,  in Problem \ref{ourproblem}. 
If the answer to Problem \ref{ourproblem} is negative,  
it might be a kind of obstruction for the existence of a consistent theory of  high dimensional Feynman diagrams (resp. that of quantization of high dimensional objects).   
If such a theory does not exist, the obstruction is a reason why it does not. 
If the answer to Problem \ref{ourproblem} is negative 
and 
if such a theory exists,  
we have to impose conditions on diffeomorphism type of high dimensional Feynman diagrams
in order to avoid the negative answer.

\smallbreak 
In \S\ref{3bdry} we prove that 
the answer to Problem \ref{ourproblem} is negative under the following condition $(\sharp\sharp)$:
Let $m\geqq3$. 
Suppose that each element of $\mathcal S$ has 
more than three connected boundary components 
(Note that each of the boundary components is a closed connected manifold).

\section{A new invariant}\label{new} 

\f We introduce a new invariant in order to discuss the $m\geqq3$ case of 
Problem \ref{ourproblem}.  

\begin{defn}\label{nu}
Let $M$ be an $m$-dimensional smooth connected compact manifold with boundary  
($m\in\N$). 
Take a handle decomposition  
\quad$A\x[0,1]\cup\mathrm{handles}$\quad
 of $M$. 
Here, we recall the following.
(See 
P.83 of \cite{Browder}, 
P.3 of \cite{Kirby}, and
\cite{Milnor, Smale}
for handle decompositions.)

\begin{enumerate}
\item 
The manifold $A$ is a closed $(m-1)$-manifold $\subset\partial M$. 
The manifold $A$ may be  $\partial M$. 
The manifold $A$ may not be  $\partial M$. 
We may have $A=\phi$.

\item 
The manifold $A$ may not be connected. 

\item 
There may be no handle (then $M=A\x[0,1]$).  
If handles are attached to $A\x[0,1]$, 
all handles are attached to $A\x\{1\}$.  
\end{enumerate}

\f Let $\mathcal H(M,A)$ denote this handle decomposition. 
An {\it ordered handle decomposition} $\mathcal H_O(M,A)$  consists of

\begin{enumerate}
\item 
a handle decomposition $\mathcal H(M,A)$ of $M$, 
 and

\item 
an order of the handles in  $\mathcal H(M,A)$: 
If we give an order to the handles, let handles be called 
 $h(\xi)$ ($\xi=1,2,3,...\delta$).  

\f The order satisfies the following: 
Let $\mu$ be a natural number $\leqq\delta$. 
Let

\noindent
$(M,A)_{\mu}=$
$A\x[0,1]$
$\begin{matrix}
\text{{\tiny $j=\mu$}}\\
\text{\large $\cup$}\\
\text{{\tiny $j=1$}}\\
\end{matrix}$
$h(j)\subset M$. 
Then this  
$A\x[0,1]$
$\begin{matrix}
\text{{\tiny $j=\mu$}}\\
\text{\large $\cup$}\\
\text{{\tiny $j=1$}}\\
\end{matrix}$
$h(j)$
is a handle decomposition of $(M,A)_{\mu}$. 
(We sometimes abbreviate $(M,A)_{\mu}$ to $M_\mu$.)
\end{enumerate}

\noindent 
For $\mu=0,$ we define $M_{\mu}=A\x[0,1]$.

\bigbreak
\f Take an ordered handle decomposition 
$\mathcal H_O(M,A)$. 
Suppose that there are ordered handles  $h(\xi)$ ($\xi=1,2,3,...\delta$).  
Let $\mu\in\{0,1,...,\delta\}$. 

\noindent
Let $\partial M_{\mu}-A\x\{0\}=E_{\mu1}\amalg...\amalg E_{\mu\xi_\mu}$, 
where each $E_{\mu i}$ is a connected closed $(m-1)$-manifold. 
Note: If $A=\phi$,   $\partial M_0-A\x\{0\}=\phi$. 
Then we suppose that  $\xi_0=1$ and $E_{01}=\phi$. 
We consider both the $A=\phi$ case and the $A\neq\phi$ case. 

\noindent
Let  
$\nu(\mathcal H_O(M,A))$ be the maximum of 
$\displaystyle\sum_{*=0}^{*=m-1}$dim$H_*(E_{\mu i};\R)$ 
for all $i, \mu$.  

\f 
Let $\nu(M,A)$ be the minimum of 
$\nu(\mathcal H_O(M,A))$ for all ordered handle decompositions 
$\mathcal H_O(M,A)$.

\f
Let $\nu(M)$ be the maximum of $\nu(M,A)$ for all $A$. 
\end{defn}

\f{\bf Note.}
By the definition,  $\nu(M)$ is an invariant of diffeomorphism type of $M$. 
If we consider  $\nu(M)$ for all smooth structures on $M$, 
we get an invariant of homeomorphism type of $M$.

\smallbreak
\f{\bf Note.} 
%
$\displaystyle\sum_{*=0}^{*=m-1}$dim$H_*(E_{\mu i};\R)$ 
is not the Euler number of $E_{\mu i}$. Their definitions are different.

\smallbreak
\f{\bf Note.} 
Suppose that we can 
attach an $m$-dimensional $k$-handle $h$ to an $m$-dimensional compact manifold $M$. 
The union $h\cup M$ is not a boundary union of $h$ and $M$ if $k<m$. 
It is a boundary union if $k=m$.


\begin{thm}\label{key}    
Let $M$ and $N$ be 
$m$-dimensional compact connected manifolds with boundaries. 
Let $M\cup_{\partial}N$ be a boundary union of $M$ and $N$. 
Then we have 
$$0\leqq\nu(M\cup_{\partial}N)\leqq \mathrm{max}\{\nu(M),\nu(N)\}.$$ 
\end{thm}

By the induction, we have a corollary.

\begin{cor}\label{keycor}    
Let $L_1,...,L_\rho$ 
 be $m$-dimensional compact connected manifolds with boundaries. 
Let 
%
%
$\begin{matrix}
\text{{\tiny $\rho$}}\\
\text{\large $\cup_{\partial}$}\\
\text{{\tiny $i=1$}}\\
\end{matrix}$
$L_i$
be 
a boundary union of  $L_1,...,L_\rho$. 
Then we have 
$$0\leqq\nu(\>
\begin{matrix}
\text{{\tiny $\rho$}}\\
\text{\large $\cup_{\partial}$}\\
\text{{\tiny $i=1$}}\\
\end{matrix}
\hskip1mm
L_i
\>)\leqq 
\mathrm{max}\{\nu(L_1),...,\nu(L_\rho)\}.$$

\end{cor}

\begin{cla}\label{claim}    
The answer to the $m\geqq3$ case of 
Problem \ref{ourproblem} 
is negative  
if 
the $\partial X=\phi$ case of 
the following Problem \ref{musttrue} has the affirmative answer.  
\end{cla}

\begin{prob}\label{musttrue}    
Let $m$ be an integer $\geqq3$. 
Suppose that there is 
an $m$-dimensional 
compact connected 
`manifold with boundary' $X$. 
Take any natural number $N$. 
Then is there 
an $m$-dimensional 
compact connected 
`manifold with boundary' $M$  
such that  $\partial M=\partial X$ 
and that  
$$\nu(M)\geqq N?$$ 
In particular, consider the $\partial X=\phi$ case. 
\end{prob}

\f{\bf{Note.}} If we do not fix the diffeomorphism type of $\partial M$, 
it is easy to prove that 
there are 
`manifolds with boundaries', $M$,   
such that $\nu(M)\geqq N.$
Because: 
Examples are 
`manifolds with boundaries', $M$,  
made from 
one 0-handle $h^0$ and $N'$ copies of 
1-handles
$h^1$, where $N'\geqq N$.

\section{Proof of Theorem \ref{key} and Claim\ref{claim}}\label{proof}

\f{\bf Proof of Theorem \ref{key}.}
By the definition of $\nu(M\cup_{\partial}N)$, 
there is a closed  $(m-1)$-manifold $P$ such that 
$$\nu(M\cup_{\partial}N)=\nu(M\cup_{\partial}N, P). \quad\cdot\cdot\cdot\cdot\cdot[[1]] $$

\noindent 
Let $A=P\cap M$. 
Let $B=P\cap N$. 
Let $C=M\cap N$. 

\bigbreak

\unitlength 0.1in
\begin{picture}(20.30,28.60)(13.10,-36.20)
%
\special{pn 8}%
\special{pa 1310 760}%
\special{pa 1310 2960}%
\special{fp}%
%
\special{pn 8}%
\special{pa 1310 2950}%
\special{pa 2010 2950}%
\special{fp}%
%
\special{pn 8}%
\special{pa 2000 2950}%
\special{pa 2000 2330}%
\special{fp}%
%
\special{pn 8}%
\special{pa 1320 770}%
\special{pa 1990 770}%
\special{fp}%
%
\special{pn 8}%
\special{pa 1980 770}%
\special{pa 2000 1930}%
\special{fp}%
%
\special{pn 8}%
\special{pa 3340 2980}%
\special{pa 3320 780}%
\special{fp}%
%
\special{pn 8}%
\special{pa 3320 790}%
\special{pa 2620 796}%
\special{fp}%
%
\special{pn 8}%
\special{pa 2630 796}%
\special{pa 2636 1416}%
\special{fp}%
%
\special{pn 8}%
\special{pa 3330 2970}%
\special{pa 2660 2976}%
\special{fp}%
%
\special{pn 8}%
\special{pa 2670 2976}%
\special{pa 2639 1816}%
\special{fp}%
%
\special{pn 8}%
\special{pa 2000 1930}%
\special{pa 2018 1904}%
\special{pa 2030 1874}%
\special{pa 2030 1842}%
\special{pa 2044 1813}%
\special{pa 2054 1783}%
\special{pa 2064 1752}%
\special{pa 2082 1726}%
\special{pa 2092 1697}%
\special{pa 2090 1664}%
\special{pa 2093 1632}%
\special{pa 2102 1601}%
\special{pa 2114 1572}%
\special{pa 2133 1546}%
\special{pa 2151 1519}%
\special{pa 2172 1495}%
\special{pa 2200 1481}%
\special{pa 2232 1472}%
\special{pa 2262 1463}%
\special{pa 2293 1452}%
\special{pa 2324 1446}%
\special{pa 2355 1443}%
\special{pa 2388 1441}%
\special{pa 2420 1438}%
\special{pa 2451 1432}%
\special{pa 2482 1424}%
\special{pa 2513 1417}%
\special{pa 2545 1411}%
\special{pa 2577 1409}%
\special{pa 2609 1410}%
\special{pa 2640 1410}%
\special{sp}%
%
\special{pn 8}%
\special{pa 2640 1830}%
\special{pa 2608 1830}%
\special{pa 2576 1830}%
\special{pa 2544 1830}%
\special{pa 2512 1830}%
\special{pa 2480 1830}%
\special{pa 2449 1831}%
\special{pa 2410 1826}%
\special{pa 2372 1824}%
\special{pa 2350 1838}%
\special{pa 2351 1870}%
\special{pa 2340 1899}%
\special{pa 2330 1929}%
\special{pa 2320 1960}%
\special{pa 2310 1990}%
\special{pa 2298 2020}%
\special{pa 2286 2050}%
\special{pa 2280 2081}%
\special{pa 2276 2112}%
\special{pa 2260 2140}%
\special{pa 2250 2170}%
\special{pa 2243 2202}%
\special{pa 2231 2231}%
\special{pa 2215 2259}%
\special{pa 2196 2285}%
\special{pa 2173 2307}%
\special{pa 2149 2331}%
\special{pa 2120 2340}%
\special{pa 2088 2340}%
\special{pa 2056 2340}%
\special{pa 2024 2340}%
\special{pa 2001 2355}%
\special{pa 2000 2360}%
\special{sp}%
\put(15.6000,-32.9000){\makebox(0,0)[lb]{$A$}}%
\put(28.9000,-33.3000){\makebox(0,0)[lb]{$B$}}%
\put(21.6000,-37.9000){\makebox(0,0)[lb]{$P=A\cup B$}}%
%
\special{pn 8}%
\special{pa 2340 1900}%
\special{pa 2030 1880}%
\special{fp}%
\put(23.3000,-21.8000){\makebox(0,0)[lb]{$C$}}%
%
\special{pn 4}%
\special{pa 2200 1920}%
\special{pa 2200 1952}%
\special{pa 2200 1984}%
\special{pa 2206 2014}%
\special{pa 2221 2043}%
\special{pa 2232 2073}%
\special{pa 2257 2094}%
\special{pa 2287 2100}%
\special{pa 2320 2100}%
\special{pa 2340 2110}%
\special{sp -0.081}%
\put(15.1000,-18.4000){\makebox(0,0)[lb]{$M$}}%
\put(28.7000,-18.8000){\makebox(0,0)[lb]{$N$}}%
%
\special{pn 4}%
\special{pa 1530 3280}%
\special{pa 1553 3303}%
\special{pa 1578 3324}%
\special{pa 1605 3342}%
\special{pa 1634 3355}%
\special{pa 1665 3363}%
\special{pa 1696 3365}%
\special{pa 1729 3363}%
\special{pa 1761 3361}%
\special{pa 1793 3359}%
\special{pa 1825 3359}%
\special{pa 1857 3360}%
\special{pa 1889 3361}%
\special{pa 1921 3363}%
\special{pa 1953 3366}%
\special{pa 1985 3371}%
\special{pa 2016 3376}%
\special{pa 2048 3380}%
\special{pa 2080 3380}%
\special{pa 2112 3378}%
\special{pa 2144 3373}%
\special{pa 2176 3368}%
\special{pa 2208 3363}%
\special{pa 2240 3358}%
\special{pa 2272 3353}%
\special{pa 2304 3350}%
\special{pa 2335 3349}%
\special{pa 2366 3351}%
\special{pa 2397 3355}%
\special{pa 2430 3360}%
\special{pa 2463 3367}%
\special{pa 2485 3390}%
\special{pa 2494 3430}%
\special{pa 2498 3468}%
\special{pa 2511 3480}%
\special{pa 2533 3455}%
\special{pa 2544 3419}%
\special{pa 2556 3393}%
\special{pa 2582 3391}%
\special{pa 2617 3402}%
\special{pa 2650 3411}%
\special{pa 2682 3417}%
\special{pa 2712 3422}%
\special{pa 2742 3427}%
\special{pa 2772 3433}%
\special{pa 2803 3441}%
\special{pa 2835 3448}%
\special{pa 2867 3452}%
\special{pa 2899 3452}%
\special{pa 2932 3448}%
\special{pa 2964 3442}%
\special{pa 2995 3433}%
\special{pa 3025 3422}%
\special{pa 3053 3410}%
\special{pa 3081 3396}%
\special{pa 3106 3381}%
\special{pa 3131 3364}%
\special{pa 3156 3345}%
\special{pa 3179 3324}%
\special{pa 3202 3301}%
\special{pa 3225 3274}%
\special{pa 3248 3246}%
\special{pa 3271 3214}%
\special{pa 3284 3183}%
\special{pa 3270 3160}%
\special{sp -0.081}%
%
\special{pn 4}%
\special{pa 2500 3500}%
\special{pa 2468 3501}%
\special{pa 2437 3492}%
\special{pa 2405 3490}%
\special{pa 2384 3509}%
\special{pa 2380 3544}%
\special{pa 2380 3576}%
\special{pa 2380 3590}%
\special{sp -0.081}%
\end{picture}%

\bigbreak

\hskip3cm Figure:  $M\cup N.$

\bigbreak

Suppose that an ordered handle decomposition 
  $\mathcal H_{O}(M, A)$ 
  gives $\nu(M, A)$. 
Hence 
$$\nu(M, A)=\nu(\mathcal H_{O}(M, A)). \quad\cdot\cdot\cdot\cdot\cdot[[2]] $$
\noindent
Let  $\mathcal H_{O}(M, A)$ 
consist of ordered handles 
$h(1),...,h(\alpha)$.

Suppose that an ordered handle decomposition 
 $\mathcal H_{O}(N, B\amalg C)$ 
 gives $\nu(N, B\amalg C)$. 
Hence 
$$\nu(N, B\amalg C)=\nu(\mathcal H_{O}(N, B\amalg C)). \quad\cdot\cdot\cdot\cdot\cdot[[3]]$$ 
\noindent
Let  $\mathcal H_{O}(N, B\amalg C)$ 
consist of ordered handles 
$k(1),...,k(\beta)$.

Let  $\mathcal H_{O}(M\cup_{\partial}N, P)$ 
be an ordered handle decomposition  
to consist of 
$l(1),...,l({\alpha+\beta})$, 
where 
the restriction of $\mathcal H_{O}(M\cup_{\partial}N, P)$ 
to 
$\left\{
\begin{array}{l}
(M, A) \\
(N, B\amalg C)
\end{array}
\right.
$
is 
$\left\{
\begin{array}{l}
\mathcal H_{O}(M, A) \\
\mathcal H_{O}(N, B\amalg C). 
\end{array}
\right.
$
That is, we have an ordered handle decomposition 

\[ 
\begin{array}{cccccccccccccc} 
M\cup_{\partial}N=&(A\amalg B)&\cup&l(1)&&\cup...&\cup&l(\alpha)&\cup&l(\alpha+1)&\cup&...\cup&l(\alpha+\beta)\\
     & ||&    &|| &   &   & & ||       & & ||  &    &    &     ||       \\
     & P &   &h(1)&   &...& & h(\alpha)& & k(1)&    &... &     k(\beta).      
\end{array}
\]

\noindent 
Here, note that 
$l(i)=h(i)(i=1,...,\alpha)$, 
and  that 
$l(i)=k({i-\alpha})(i=\alpha+1,...,\alpha+\beta)$.

\bigbreak
Recall $E_{\mu i}$ in Definition \ref{nu}.  
Take  $E_{\mu i}$ for this 
$\begin{cases}
\mathcal H_{O}(M\cup_{\partial}N, P)\\ 
\mathcal H_{O}(M, A) \\
\mathcal H_{O}(N, B\amalg C), 
\end{cases}$
call it 
$\begin{cases}
E^{M\cup_{\partial}N}_{\sharp\star} \\ 
E^{M}_{\clubsuit\diamondsuit}  \\
E^{N}_{\heartsuit\spadesuit}.   
\end{cases}$

\noindent
Then 
$\{E^{M\cup_{\partial}N}_{\sharp\star}| \text{$\sharp$ and $\star$ take all values}\}=$\newline
 $\{E^{M}_{\clubsuit\diamondsuit}| \text{$\clubsuit$ and $\diamondsuit$ take all values}\}
\cup\{E^{N}_{\heartsuit\spadesuit}| \text{$\heartsuit$ and $\spadesuit$ take all values}\}$.

\bigbreak
\noindent  
{\bf Note.}  
Furthermore we have the following.

\vskip1mm\noindent
If $0\leqq\mu\leqq\alpha$,   
$\{E^{M\cup_{\partial}N}_{\mu\star}| \text{$\star$ takes all values}\}=$ \newline
$\{E^{M}_{\mu\diamondsuit}| \text{$\diamondsuit$ takes all values}\}\cup$ 
$\{B_\flat | \text{$\flat$ takes all values}\}$. 
Here 
$B=B_1\amalg...\amalg B_\zeta$, 
where each $B_\flat$ is a closed connected manifold. 
$B_\flat\in$$\{E^{N}_{0\clubsuit}| \text{$\clubsuit$ takes all values}\}$ holds.

\vskip1mm\noindent
If $\alpha+1\leqq\mu\leqq\alpha+\beta$,   
$\{E^{M\cup_{\partial}N}_{\mu\star}| \text{$\star$ takes all values}\}=$ \newline
$\{E^{N}_{{\mu-\alpha}\hskip1mm\clubsuit}| \text{$\clubsuit$ takes all values}\}\cup$ 
$\{D_\natural |\text{$\natural$ takes all values}\}$. 
Here, $\partial M- A - C=$  
$D_1\amalg...\amalg D_\varepsilon$, 
where each $D_\natural$ is a closed connected manifold.  
$D_\natural\in$ $\{E^{M}_{\alpha\diamondsuit}| \text{$\diamondsuit$ takes all values}\}$ holds.



\bigbreak
Suppose 
$\nu(\mathcal H_{O}(M\cup_{\partial}N, P))=$ 
$\displaystyle\sum_{*=0}^{*=m-1}$dim$H_*(E^{M\cup_{\partial}N}_{\mu i};\R)$ 
for an integer $\mu$ and an integer $i$.   
Then 
 $E^{M\cup_{\partial}N}_{\mu i}$ is \newline 
$E^{M}_{\sigma j}$ for an integer $\sigma$ and an integer $j$ $\quad\cdot\cdot\cdot\cdot\cdot$[[I]]
\newline 
or  \newline  
$E^{N}_{\tau k}$ for an integer $\tau$ and an integer $k$.   $\quad\cdot\cdot\cdot\cdot\cdot$[[II]]

\noindent
Suppose [[I]] holds. Then 
$\nu(\mathcal H_{O}(M, A))=
\displaystyle\sum_{*=0}^{*=m-1}$dim$H_*(E^{M}_{\sigma j};\R)$ 
for the integer $\sigma$ and the integer $j$.   Hence 
$$\nu(\mathcal H_{O}(M\cup_{\partial}N, P))=\nu(\mathcal H_{O}(M, A)).  \quad\cdot\cdot\cdot\cdot\cdot[[4]]$$ 

\noindent
Suppose [[II]] holds. Then 
$\nu(\mathcal H_{O}(N, B\amalg C))=
\displaystyle\sum_{*=0}^{*=m-1}$dim$H_*(E^{N}_{\tau k};\R)$ 
for the integer $\tau$ and the integer $k$.   Hence 
$$\nu(\mathcal H_{O}(M\cup_{\partial}N, P))=\nu(\mathcal H_{O}(N, B\amalg C)).  \quad\cdot\cdot\cdot\cdot\cdot[[5]]$$

\smallbreak
\f 
Since 
[[4]] or
[[5]]  holds, 
$$\nu(\mathcal H_{O}(M\cup_{\partial}N, P))
\leqq 
\mathrm{max}  \{\nu(\mathcal H_{O}(M, A)), 
\nu(\mathcal H_{O}(N, B\amalg C))\}.  \quad\cdot\cdot\cdot\cdot\cdot[[6]] $$

\bigbreak

By [[2]], [[3]], and [[6]], 
$$\nu(\mathcal H_{O}(M\cup_{\partial}N, P))
\leqq{\rm max}\{\nu(M, A), \nu(N, B\amalg C)\}.  \quad\cdot\cdot\cdot\cdot\cdot[[7]]$$

\f 
By the definition of $\nu(M\cup_{\partial}N, P),$  
$$\nu(M\cup_{\partial}N, P)
\leqq\nu(\mathcal H_{O}(M\cup_{\partial}N, P)).  \quad\cdot\cdot\cdot\cdot\cdot[[8]]$$

\noindent 
By [[7]] and [[8]], 
$$\nu(M\cup_{\partial}N, P)\leqq 
{\rm max}\{\nu(M, A), \nu(N, B\amalg C)\}. 
\quad\cdot\cdot\cdot\cdot\cdot[[9]]$$

\f
By [[1]] and [[9]]
$$\nu(M\cup_{\partial}N)
\leqq
{\rm max}\{\nu(M, A), \nu(N, B\amalg C)\}. \quad\cdot\cdot\cdot\cdot\cdot[[10]]$$ 

\f 
By the definition of $\nu(M)$ and $\nu(N)$,  we have 
$$\nu(M,A)
\leqq
\nu(M)  
\quad{\rm and}\quad  
\nu(N, B\amalg C) 
\leqq 
\nu(N). \quad\cdot\cdot\cdot\cdot\cdot[[11]]$$ 

\f 
By [[10]] and [[11]], 
$$\nu(M\cup_{\partial}N)
\leqq  
{\rm max}\{\nu(M), \nu(N)\}.$$

\smallbreak

By the definition of $\nu(\quad)$, 
 $0\leqq\nu(M\cup_{\partial}N)$. 
This completes the proof. 
\qed

\bigbreak
\f{\bf Proof of Claim \ref{claim}.}
We suppose the following assumption and we deduce a contradiction.

Assumption: we have the affirmative answer to Problem \ref{ourproblem}.

By the above assumption there is a finite set 
$\mathcal S=\{S_1,...,S_s\}$ as in Problem \ref{ourproblem}.

\f 
By Corollary \ref{keycor},   
for any connected closed  $m$-manifold $L$, 
$\nu(L)\leqq$max$\{\nu(S_1),...,\nu(S_s)\}$.

If 
the $\partial X=\phi$ case of Problem \ref{musttrue} has the affirmative answer, 
then there is 
a connected closed  $m$-manifold $M$  
such that 
$\nu(M)>$max$\{\nu(S_1),...,\nu(S_s)\}$. 

We arrived at a  contradiction. Hence Claim \ref{claim} is true. 
\qed

\section{Some results on our new invariant}\label{someresults}

\f 
Let $M\neq\phi.$ 
Let $m\geqq2.$ 
Let $M$ be a smooth closed oriented connected $m$-manifold.
By the definition,  $\nu(M)\geqq2$.

We prove the following.

\smallbreak
\begin{thm}\label{sphere}    
Let $m\geqq2.$ 
Let $S^m$ be diffeomorphic to the standard sphere. 
Then $\nu(S^m)=2$.
\end{thm}

\f{\bf Proof of Theorem \ref{sphere}.}
There is an ordered handle decomposition 
$\mathcal H_O=h(1)\cup h(2)$ 
such that 
$h(1)=h^0$, 
$h(1)=h^m$.
Then $\nu(\mathcal H_O)=2$. 
Hence 
$\nu(S^m)\leqq2$. 
By the definition, $\nu(S^m)\geqq2$. 
Hence $\nu(S^m)=2$.

\f{\bf{Note.}}
Furthermore we have the following: 
If $M$ has a handle decomposition $h^0\cup h^m$, 
then $\nu(\mathcal H)=2$.

\bigbreak

It is natural to ask the following:
Suppose $M$ is a closed connected oriented manifold. 
Then does $\nu(M)=2$ imply that $M$ is PL homeomorphic to $S^m$?

\f
We have the following theorem as an answer to this question.

\begin{thm}\label{Chicago1}  
Let $n\geqq2.$
We have 
$\nu(\>
\begin{matrix}
\text{{\tiny $\star$}}\\
\text{\large $\sharp$}\\
\text{{\tiny $$}}\\
\end{matrix}
(S^1\times S^{n-1})
\>)=2$, 
where $\star$ is any natural number
and where 
$\begin{matrix}
\text{{\tiny $1$}}\\
\text{\large $\sharp$}\\
\text{{\tiny $$}}\\
\end{matrix}
(S^1\times S^{n-1})
%
=S^1\times S^{n-1}$.  
\end{thm}

\noindent
{\bf Proof of Theorem \ref{Chicago1}.}  
There is an ordered handle decomposition 

\noindent
$h^0\cup h^{n-1}\cup h^1\cup h^{n-1}\cup h^1\cup.....\cup h^{n-1}\cup h^1\cup h^n$,  
where this order is the order of this ordered handle decomposition. 
Recall $E_{\mu i}$ in Definition \ref{nu}. 
We can suppose that each $E_{\mu i}$ is a sphere. 
\qed

\bigbreak 
We have more results on the $\nu$ invariant. 

\begin{thm}\label{notequal}
In Theorem \ref{key}, there is a pair of `manifolds with boundaries' $M$, $N$  
such that 
$$\nu(M\cup_{\partial}N)\neq \mathrm{max}\{\nu(M),\nu(N)\}.$$ 
\end{thm}

\f{\bf Proof of Theorem \ref{notequal}.} Let $M\cong N$. 
Let $M$ be the 2-dimensional solid torus $S^1\x D^2$. 
Note $\partial M$ is the torus $T^2$. 
The closed manifold $A$ in Definition \ref{nu} for $S^1\x D^2$ is $T^2$ or $\phi$.  
Consider all 
$\mathcal H_O(S^1\x D^2 ,T^2)$
 and $\mathcal H_O(S^1\x D^2 ,\phi)$.  
Then all these $\mathcal H_O$ have 
$\partial M_\mu=T^2$ for an integer $\mu$.  
Hence $\nu(\mathcal H)\geqq4$. 

There is a handle decomposition  $h^0\cup h^1$ such that 
$\nu(\mathcal H_O(S^1\x D^2, \phi))=4$. 

\f  There is a handle decomposition   
$T^2\x[0,1]\cup h^2\cup h^3$ such that  
$\nu(\mathcal H_O(S^1\x D^2, T^2))=4$. Hence $\nu(M)=\nu(N)=4$.

Note that there is a boundary union $S^3=M\cup_\partial N$.
By Theorem \ref{sphere},  $\nu(S^3)=2$. 

\f Hence 
$\nu(M\cup_{\partial}N)=2<4=\mathrm{max}\{\nu(M),\nu(N)\}.$ 

\f Hence  $\nu(M\cup_{\partial}N)\neq \mathrm{max}\{\nu(M),\nu(N)\}.$ 
\qed

\bigbreak
It is natural to ask the following question. 
Is there a closed connected oriented manifold $M$ with $\nu(M)>2$?

We have the following theorems as answers to this question.

\begin{thm}\label{ChicagoA}
Let $M$ be a connected closed oriented 3-manifold. 
Suppose that $\pi_1(M)$ is not the trivial group or 
any of the free groups $\Z*...*\Z$. 
Then   $\nu(M)>2$. 
\end{thm}

\begin{thm}\label{ChicagoB}      
Let $M$ be one of the 3-dimensional lens spaces $L(p,q)$, 
where $L(p,q)$ is not $S^3$ or $S^1\times S^2$ as usual. 
Then   $\nu(M)=4$. 
\end{thm}

\noindent{\bf Note to Theorem \ref{ChicagoB}.}  $\pi_1(M)=\Z_p$.   

\smallbreak

\noindent{\bf Proof of Theorem \ref{ChicagoA}.}  
We suppose the following assumption and we deduce a contradiction. 

Assumption.   $\nu(M)\leqq 2$.

By the definition of the $\nu$-invariant, $\nu(M)\geqq2$. 
Hence $\nu(M)\leqq2$ and $\nu(M)\geqq2$. 
Hence $\nu(M)=2$. 

Take an ordered handle decomposition 
$h(1)\cup h(2)\cup...\cup h(\mu)\cup...$
which determines $\nu(M)=2$.  
Take $M_\mu$ and $E_{\mu i}$ in Definition \ref{nu}. 
Let $\mu\geqq1$ because $M_\mu=\phi$ if $\mu=0$. 
Since $M$ is a connected closed oriented 3-manifold, 
each $E_{\mu i}$ is a connected closed orientable surface. 
Hence  $\displaystyle\sum_{j=0}^{j=2}{\text{dim}} H_j(E_{\mu i}:\R)\geqq2$ 
for each $i$. 
By $\nu(M)=2$ 
we have \newline
$\displaystyle\sum_{j=0}^{j=2}{\text{dim}} H_j(E_{\mu i}:\R)\leqq2$ 
for each $i$. 
Hence 
$\displaystyle\sum_{j=0}^{j=2}{\text{dim}} H_j(E_{\mu i}:\R)=2$ 
for each $i$. 
Hence each $E_{\mu i}$ is a 2-sphere 
for each $i$. 

\bigbreak
Let $M_\mu=M_{\mu1}\amalg...\amalg M_{\mu\psi}$, 
where $\amalg$ denotes a disjoint union and 
where each $M_{\mu\star}$ is a connected compact manifold with boundary. 

Let $\Pi(M_\mu)$ be a `finite sequence of groups'
$(\pi_1(M_{\mu1}),...,\pi_1(M_{\mu\psi}))$.

\bigbreak

Suppose that $h(\mu+1)$ is a 3-handle. 
Note that $M_{\mu+1}=M_\mu\cup h(\mu+1)$. 
By Van Kampen's theorem, 
$\Pi(M_{\mu+1})=\Pi(M_{\mu})$.

\bigbreak
Suppose that $h(\mu+1)$ is a 2-handle. 
Then the core of the attaching part of $h(\mu+1)$  
is a circle. 
Since the attaching part is connected, 
the attaching part  is included  in one of the 2-spheres $E_{\mu i}$, 
call it $E_{\mu \iota}$. 
Hence it holds that    

\smallbreak\hskip1cm
the attaching part of $h(\mu+1)\subset$ the 2-sphere $E_{\mu \iota}$ $\subset M_\mu$. 
\smallbreak
\noindent 
Hence the attaching part is contractible in $M_\mu$. 
By Van Kampen's theorem, 
\newline$\Pi(M_{\mu+1})=\Pi(M_{\mu})$.

\bigbreak

Suppose that $h(\mu+1)$ is a 1-handle. 
One of the following two conditions holds if we change the suffix $*$ of $G_*$. 

\noindent
(1)
$\Pi(M_{\mu})=(G_1,..., G_{\psi-1}, G_\psi)$ and 
$\Pi(M_{\mu+1})=(G_1,..., G_{\psi-1}, G_\psi*\Z)$.  

\noindent
(2)
$\Pi(M_{\mu})=(G_1,..., G_{\psi-2}, G_{\psi-1}, G_\psi)$ and 
$\Pi(M_{\mu+1})=(G_1,..., G_{\psi-2}, G_{\psi-1}* G_\psi)$.  

\bigbreak
Suppose that $h(\mu+1)$ is a 0-handle. 
The following condition holds if we change the suffix $*$ of $G_*$. 
$\Pi(M_{\mu})=(G_1,...,G_\psi)$ and 
$\Pi(M_{\mu+1})=(G_1,...,G_\psi,1)$. 
Here, 1 denotes the trivial group.

\bigbreak
Note that $M_0=\phi$ and that $M_1=h(1)$ is a 0-handle. 
Hence $\pi_1(M)$ is the trivial group or 
one of the free groups $\Z*...*\Z$.   
This is a contradiction. 

This completes the proof. 
\qed
\smallbreak

\noindent{\bf Proof of Theorem \ref{ChicagoB}.}   
By Theorem \ref{ChicagoA}   
and $\pi_1M=\Z_p$, 
we have $\nu(M)>2$. 

As in the Proof of Theorem \ref{ChicagoA},  
each $E_{\mu i}$ is a connected closed orientable surface.  
Hence  
$\displaystyle\sum_{j=0}^{j=2}{\text{dim}} H_j(E_{\mu i}:\R)$ is an even number 
for each $\mu\geqq1$ and each $i$. 
Hence $\nu(M)$ is an even number. 
By $\nu(M)>2$, we have $\nu(M)\geqq4$.

There is an ordered handle decomposition $h^0\cup h^1\cup h^2\cup h^3$ of $M$ 
such that 
$\partial M_1=S^2$, 
$\partial M_2=T^2$, 
$\partial M_3=S^2$, 
$\partial M_4=\phi$. 
Hence $\nu(M)\leqq4$. 

Since $\nu(M)\geqq4$ and $\nu(M)\leqq4$, 
we have $\nu(M)=4$. 

This completes the proof. 
\qed

\smallbreak\noindent
Alternative proof of $\nu(L(p,q))\leqq4$. 
Note that $L(p,q)$ is a boundary union of two solid torus. 
 In the proof of Theorem \ref{notequal},   
we prove that 
the $\nu$ invariant of the solid torus is four. 
By Theorem \ref{key}, 
$\nu(L(p,q))\leqq4$.

\section{The solution to Problem 1.2 in a special case}\label{3bdry}

As a partial solution to Problem \ref{ourproblem}, 
we prove that 
the answer to Problem \ref{ourproblem} is negative under the condition $(\sharp\sharp)$ 
in the last paragraph of \S\ref{introduction}.

In this case we obtain the result without using the $\nu$-invariant.

\bigbreak 
We suppose that the following assumption is true,  
and deduce a contradiction. 

\smallbreak\f {\bf Assumption.} 
Then the answer to Problem \ref{ourproblem} is affirmative 
under the condition $(\sharp\sharp)$ 
in the last paragraph of \S\ref{introduction}.

\smallbreak Let $W$ be an $m$-dimensional arbitrary compact connected manifold 
with boundary.  
Suppose that $\partial W$ has $z$ connected components.  Note that we fix $z$. 
We prove both the $z=0$ case and the $z\geqq1$ case. Note that  
we suppose that $z=0$  in Problem \ref{ourproblem} and that 
$z\geqq1$ in the paragraph right under Problem \ref{ourproblem}

By the assumption we can divide $W$ into pieces $W_i\in \mathcal S$ 
and can regard $W=W_1{\cup_\partial}...{\cup_\partial} W_w$.

\smallbreak
Let 
$W_i \cap W_{i'}$ have  $\rho$ connected components. 
Hence  $\rho\geqq{\frac{3w-z}{2}}$.

\smallbreak

Consider the Meyer-Vietoris exact sequence:

\f
$H_j( \amalg_{i,i'}\{W_i \cap W_{i'}\};\R)
\to H_j(\amalg_{i=1}^{i=w} W_i;\R)\to H_j(W;\R)$. 
Here, $\amalg_{i,i'}$ denotes the disjoint union of $W_i \cap W_{i'}$ for all $(i,i')$.
Consider 

\f 
$H_1(W;\R)\to H_0(\amalg_{i,i'}\{W_i\cap W_{i'}\};\R)\to
H_0(\amalg W_i;\R)\to H_0(W;\R)\to0$. 
Note 

\f
$H_0(\amalg W_i;\R)\cong\R^w$ and $H_0(W;\R)\cong\R$. 
We suppose $H_1(W;\R)\cong\R^l$. 
Then we have the exact sequence: 

\f
$\R^l\to\R^\rho\to\R^w\to\R\to0.$   
Hence $l\geqq\rho-w+1$. 
Hence $l\geqq\frac{w-z+2}{2}$. 
Hence $(2l+z-2)\geqq w$.

\smallbreak
We define an invariant $h(\quad)$. 
Let $X$ be a compact manifold with boundary. 
Take a handle decomposition of $X$. 
Consider the numbers of handles in the handle decompositions of $(X,A)$, 
where $A$ is one as defined in Definition \ref{nu}. 
Let $h(X,A)$ be the minimum of such the numbers.  
Let $h(X)$ be the maximum of $h(X,A)$ for all $A$.

Note that 
$\mathcal S$ is a finite set $\{M_1,...,M_\mu\}$. 
Suppose that $M'$ is one of 
the `manifolds with boundaries' $M_i$ and 
that $h(M')\geqq h(M_i)$ for any $i$.  
Then we have $w\times h(M')\geqq h(W)$. 
Hence $(2l+z-2)\times h(M')\geqq h(W)$. 
Note that the left side is constant. 

\smallbreak

For any natural number $N$, 
there are countably infinitely many compact oriented connected 
$m$-manifolds $W'$ with boundaries 
such that $\partial W'=\partial W$, 

\f
that 
$H_1(W;\R)\cong\R^l$, 
and that $h(W)\geqq N$. 
Because: There is an $n$-dimensional closed manifold $P$ such that 
$H_1(P;\R)\cong\R^l$. 
There is an $n$-dimensional rational homology sphere $Q$ 
which is not an integral homology sphere. 
Make a connected sum which is made from one copy of $P$ and $q$ copies of $Q$
($q\in\N\cup\{0\}$).

We arrived at a contradiction.  This completes the proof. \qed


\bigbreak
Furthermore \cite{Ogasa} pointed out  the following. 

\f
(1) 
There is an  
$n$-dimensional connected Feynman diagram with three outlines
(a compact connected $n$-manifold with boundary whose boundary has three connected components)   
with the following properties. 
Two copies of the diagram is made into 
countably infinitely many kinds of 
diagrams 
with four outlines. 

\f(The idea of the proof: Let the diagram be 
$\{($the solid torus$)-$two open 3-balls$\}$.  
Use the fact that all 3-dimensional lens spaces, $S^3$, and $S^1\x S^2$ 
are  made from two solid torus.)

\smallbreak
\f
(2) 
There is an infinite set  $\mathcal S$ with the following properties. 

\noindent
(i) All $m$-dimensional Feynman diagrams 
(compact $m$-manifolds with boundaries) 
are boundary unions of finite elements of $\mathcal S$.  

\noindent
(ii) Each element of $\mathcal S$ is 
 what is made by attaching an $m$-dimensional handle to 

\noindent
(an $(m-1)$-dimensional connected closed manifold)$\x[0,1]$. 
Note it has one, two or three connected boundary components. 
%
(The idea of the proof: Use handle decompositions.)

\section{Discussion}\label{discussion}

\f Take a group 
$G=\{g_1,...,g_N|
\newline 
 g_1\cdot g_2...\cdot g_{N-1}\cdot g_N\cdot g_2^{-1}...\cdot g_{N}^{-1}=1, 
 \newline 
 g_2\cdot g_3...\cdot g_{N}\cdot g_1\cdot g_3^{-1}...\cdot g_{N}^{-1}\cdot g_1^{-1}=1,..., 
\newline 
 g_N\cdot g_1...\cdot g_{N-2}\cdot g_{N-1}
\cdot g_1^{-1}...\cdot g_{N-1}^{-1}=1. \}$

In the  $m\geqq4$ case, 
we can  make a compact connected oriented manifold $Z$ with boundary such that 

\f(1) $\pi_1(Z)=G$. 

\f(2) $Z$ is made of one 0-handle, 
$N$ copies of 1-handles, and 
$N$ copies of 2-handles.
(Each of the generators $g_*$ corresponds to each of the 1-handles.  
Each of the $N$ relations corresponds to each of  the 2-handles.)

\f Take the double of Z. Call it $W$. Note $\pi_1(W)=G.$


\noindent
Thus we submit the following problem. 

\begin{prob}\label{hope}
Do you prove $\nu(W)\geqq N$?  
\end{prob}

\f
If the answer to Problem \ref{hope} is affirmative, 
then the answer to Problem \ref{musttrue} is affirmative
(in the closed manifold case, which would be extended in all cases). 

\f
By using a manifold with boundary whose fundamental group is so complicated as above, 
we may solve Problem \ref{ourproblem}, \ref{musttrue}.

We explain the above strategy more.

\bigbreak

Let $M$ be a finite dimensional smooth compact connected manifold with boundary. 

\noindent
Take any handle decomposition $\mathcal H$, 

\noindent 
$M=h^0\cup h^1_1...\cup h^1_\xi$ $\cup h^2$'s.....

\noindent 
Suppose that there is only one 0-handle $h^0$ 
in this handle decomposition $\mathcal H$. 

\noindent 
Suppose that the 1-handles $h^1_1...h^1_\xi$ are 
all 1-handles of this handle decomposition $\mathcal H$. 

\noindent 
Different 1-handles 
$h^1_{\alpha_1},...,h^1_{\alpha_k}\in\{h^1_1,...,h^1_\xi\}$ 
are called {\it brother-handles} 
if  $h^1_{\alpha_1},...,h^1_{\alpha_k}$ satisfy the following. 
If a 2-handle of the handle decomposition $\mathcal H$
is attached to one of  $h^1_{\alpha_1},...,h^1_{\alpha_k}$, 
then the 2-handle is attached to all  of  $h^1_{\alpha_1},...,h^1_{\alpha_k}$.  

Take each set of brother-handles of $\mathcal H$. 
Take the order of each set. 
Take the maximum of the orders. 
Call it $b(\mathcal H)$.

Consider each handle decomposition $\mathcal E$ of $M$. 
Take $b(\mathcal E)$ for each handle decomposition $\mathcal E$.
Take the minimum of all $b(\mathcal E)$.
Call it $b(M)$.

We submit the following problem.

\noindent
\begin{prob}\label{tsuika1}
Let $m$ be an integer $\geqq3$. 
Suppose that there is 
an $m$-dimensional 
compact connected manifold $X$ with boundary. 
Take any natural number $n$. 
Then is there 
an $m$-dimensional 
compact connected manifold $M$ with boundary 
such that  $\partial M=\partial X$ 
and that  
$$b(M)\geqq n?$$ 
In particular, consider the $\partial X=\phi$ case. 
\end{prob}

If Problem \ref{tsuika1} has the affirmative answer, 
then Problem \ref{musttrue} may have the affirmative answer.

The author could prove that 
$b(S^n)=0$,  
$b(\R P^n)=1$,  
$b(\sharp^m\R P^n)=1$, 
$b(T^2)=2$ and some other cases. 
The author guesses  that 
it might not be difficult to calculate  the $b(\quad)$ 
in the following Problem \ref{tsuika2} 
or to evaluate it from the lower side.

\noindent
\begin{prob}\label{tsuika2} 
For the closed manifold $W$ in Problem \ref{hope}, 
is $b(W)\geqq N$? 
\end{prob}

The author would think that these problems could be solved by group theoretic ways. 

The author tries to interpret Problem\ref{tsuika1} 
in terminology of  group theory.

Let $G$ be a finitely generated group. 
Let $<g_1,...,g_\xi|r_1,...,r_\zeta>$ be a presentation $\mathcal P$ of $G$. 
Different generators 
$g_{\alpha_1},...,g_{\alpha_k}\in\{g_1,...,g_\xi\}$ 
are called {\it brother-generators} 
if  $g_{\alpha_1},...,g_{\alpha_k}$ satisfy the following. 
If a relation $r_*$ of the presentation $\mathcal P$ 
includes one of  $g_{\alpha_1},...,g_{\alpha_k}$, 
then the relation $r_*$ includes all  of  $g_{\alpha_1},...,g_{\alpha_k}$.  
Here, if we say that  $r_*$ includes $g_\rho$, then 
it means that $r_*$ includes $g_\rho$ or $(g_\rho)^{-1}$. 

Take each set of brother-generators of $\mathcal P$. 
Take the order of each set. 
Take the maximum of the orders. 
Call it $b(\mathcal P)$.

Consider each presentation $\mathcal E$ of $G$. 
Take $b(\mathcal E)$ for each presentation $\mathcal E$.
Take the minimum of all $b(\mathcal E)$.
Call it $b(G)$.

We submit the following problems.

\noindent
\begin{prob}\label{tsuika3}
Let $n$ be any natural number. 
For a group $G$, 
is $b(G)$ greater than $n$? 
\end{prob}

The following problem may be connected with Problem \ref{tsuika2}. 

\noindent
\begin{prob}\label{tsuika4}
For the closed manifold $W$ in Problem \ref{hope}, 
do we have $b(\pi_1(W))\geqq N$? 
\end{prob}

We submit one more problem which seems easier than other problems in this paper.
\noindent
\begin{prob}\label{easy}
Is the $\nu$ invariant of the Poincar\'e sphere six? 
\end{prob}

Use $Z_p$ coefficient homology groups instead in the definition of $\nu$. 
Use the order of Tor$H_*(\quad;\Z)$ instead in the definition of $\nu$. 
Can we solve Problem \ref{ourproblem}?

\smallbreak 
Calculate $\nu$, $b$ of the knot complement. 
In particular, in the case of 1-dimensional prime knots. 
In this case, what kind of connection with the Heegaard genus?

\smallbreak 
If we replace 
$\displaystyle\sum_{*=0}^{*=m-1}$dim$H_*(E_{\mu i};\R)$ in Definition \ref{nu}
with 
dim$H_*(E_{\mu i};\R)$  for the fixed integer $*$  
or 
a nonnegative real number valued topological invariant  
(resp. diffeomorphism type invariant) of $E_{\mu i}$,    
we obtain another invariant instead of $\nu$.   
It satisfies Theorem \ref{key}.

\smallbreak 
\noindent
An example made by using such a nonnegative real number valued topological invariant:
If we replace 
$\displaystyle\sum_{*=0}^{*=m-1}$dim$H_*(E_{\mu i};\R)$  in the m=4 case in Definition \ref{nu} 
with 
 the absolute value of a quantum invariant $\tau_\star(E_{\mu i})$,   
we obtain another invariant instead of $\nu$. 
It satisfies Theorem \ref{key}.

\smallbreak 
This paper is based on the author's preprints \cite{Ogasa}.


\vskip3mm
\f 
Computer Science, Meijigakuin University, Yokohama, Kanagawa, 244-8539, Japan 

\f pqr100pqr100@yahoo.co.jp, ogasa@mai1.meijigkauin.ac.jp

\end{document}